\documentclass[12pt,reqno]{amsart}
\usepackage{amsmath,amssymb,amsthm}
\usepackage{geometry}
\newtheorem{theorem}{Theorem}[section]
\newtheorem{lemma}[theorem]{Lemma}
\newtheorem{corollary}[theorem]{Corollary}

\newtheorem*{thmA}{Theorem A}
\newtheorem*{thmB}{Theorem B}

\title{On finite simple groups possessing a conjugacy class of size $3p^n$}

     \author[Beltr\'an and Shao]{Antonio Beltr\'an\\
     Departamento de Matem\'aticas\\
      Universitat Jaume I \\
     12071 Castell\'on\\
      Spain\\ 
      \\Changguo Shao\\
      College of science\\
NANJING UNIVERSITY OF POSTS AND TELECOMMUNICATIONS\\
NANJING 210023 YADONG
        \\
      China}

 \thanks{Antonio Beltr\'an: abeltran@uji.es,
 Changguo Shao: shaoguozi@163.com}



\date{}

\begin{document}

\begin{abstract}
We  classify all non-abelian finite simple groups that  possess conjugacy classes of size $3p^n$ for some prime $p$ and non-negative integer $n$. As an application, we study the structure of the normal closure $\langle g^G\rangle$ of an element $g$ in an arbitrary finite group $G$, provided that the conjugacy class $g^G$ has size $3p^n$. Specifically, we provide  conditions for $\langle g^{G}\rangle$ to be solvable and completely determine its quotient modulo the solvable radical in the non-solvable case.
\end{abstract}

\maketitle

\section{Introduction}

The study of the influence of local numerical invariants of conjugacy classes on the structure of finite groups has a rich history, dating back to the seminal result of Burnside, who proved that any group possessing a conjugacy class of prime-power size cannot be non-abelian simple.  Over the decades, a considerable number of authors have generalized and expanded this line of research, investigating how restrictions on the prime factorization of conjugacy class sizes restrict solvability, nilpotency, or structural properties (see, for instance, \cite{Huppert} or the survey by Camina and Camina \cite{Camina}). A major milestone in this context was achieved by Kazarin \cite{Kaz}, who proved, without resorting to the Classification of Finite Simple Groups (CFSG), that if an element $g$ of a finite group $G$ has a conjugacy class of prime-power size, then $g$ belongs to the solvable radical of $G$,  thereby implying that $\langle g^G \rangle$ is solvable. Building upon Kazarin's work, the authors recently investigated in \cite{BS} the arithmetic condition where a conjugacy class has size $2p^n$ for some prime $p$ and integer $n \ge 0$. We established that no non-abelian finite simple group admits  such a conjugacy class. This resulted in  a description of the  structure of  $\langle g^G \rangle$ for an element $g$ in a finite group $G$ with class size  $2p^n$ when $\langle g^G \rangle$ is non-solvable.

\medskip
In the present paper, we go further by investigating the structure of finite groups possessing a conjugacy class of size $3p^n$, where $p$ is a prime and $n \ge 0$. Extending the arithmetic bound from $2p^n$ to $3p^n$ leads to a more delicate situation: unlike the $2p^n$ case, non-abelian finite simple groups possessing conjugacy classes of size $3p^n$ do exist. However, such instances are rare. Specifically, these simple groups and their corresponding conjugacy classes are completely classified in our first main result.

\begin{thmA}\label{Theorem A}
Let $G$ be a  finite non-abelian  simple group that possesses a conjugacy class $g^G$ of size $3p^n$ for some prime $p$ and integer $n \ge 0$. Then $G  \cong A_5$ or $G \cong {\rm PSL}(2,7)$. Furthermore, the explicit cases are:
\begin{enumerate}
\item[(1)] $G \cong A_5$, $|g^G| = 15 = 3 \cdot 5$, and $o(g) = 2$ (so $p = 5$);
    \item[(2)] $G \cong A_5$, $|g^G| = 12 = 3 \cdot 2^2$, and $o(g) = 5$ (so $p = 2$);
    \item[(3)] $G \cong {\rm PSL}(2, 7)$, $|g^G| = 21 = 3 \cdot 7$, and $o(g) = 2$ (so $p = 7$);
    \item[(4)] $G \cong {\rm  PSL}(2, 7)$, $|g^G| = 24 = 3 \cdot 2^3$, and $o(g) = 7$ (so $p = 2$).
\end{enumerate}
\end{thmA}

The proof of Theorem  A relies heavily on  CFSG and makes essential use of the classification   of subgroups with indices  divisible by at most two primes, established  by Li and Li (\cite{LL}). We then apply Theorem A to general finite groups  to study the structure of $\langle g^G \rangle$ when the size of $g^G$ is $3p^n$. When $\langle g^G \rangle$ fails to be solvable, we determine the structure for the factor group $\langle g^G \rangle / S(\langle g^G \rangle)$, where $S(\langle g^G \rangle)$ is the solvable radical, demonstrating that this group must be isomorphic to one of a short, explicit list of almost simple groups with socle $A_5$, $A_6$, $A_7$, or ${\rm PSL}(2,7)$. 

\begin{thmB}
Let $G$ be a finite group having a conjugacy class $g^G$ of size $3p^n$, where $p$ is prime and $n \ge 1$, such that $\langle g^G \rangle$ is non-solvable. Let $\overline{g}$ denote the image of $g$ in $\langle g^G\rangle/S(\langle g^G\rangle)$. Then  $\langle g^G\rangle/S(\langle g^G\rangle)$ is isomorphic to one of the following  groups:
\begin{enumerate}
    \item $A_5$, where $\overline{g}$ is an element of order $2$ or $5$.
    \item $S_6$, where $\overline{g}$ is an odd permutation of order $2$ in $S_6$ (either a transposition or a product of three disjoint transpositions).
    \item $S_7$, where $\overline{g}$ is a transposition in $S_7$.
    \item $\mathrm{PSL}(2, 7)$, where $\overline{g}$ is an element of order $7$, or an element of order $2$ lying in the center of a Sylow $2$-subgroup of $\mathrm{PSL}(2, 7)$.
   \end{enumerate}
\end{thmB}

This illustrates how arithmetic restrictions on the length of a single conjugacy class can impose strong constraints on the global structure of a finite group.

\section{Proof of Theorem A}

The proof of Theorem A is divided into three parts according to the three families of simple groups in the CFSG (Theorems \ref{thm1.1}, \ref{thm1.2} and \ref{thm1.3}).

\begin{theorem}\label{thm1.1}
The alternating group $A_m$ with $m\geq 5$ has no conjugacy class of size $3p^n$ for any prime $p$ and integer $n \ge 0$, except in the case 
$A_5$, with
\begin{itemize}
    \item[(a)]
$p=5, n=1$,  and class size $15 = 3 \cdot 5$, given by elements of order $2$.
\item[(b)] $p=2, n=2$  and class size $12 = 3 \cdot 2^2$, given by elements of order $5$.
   \end{itemize} 
\end{theorem}

\begin{proof}
A direct verification, employing the Atlas \cite{Con} or GAP \cite{GAP4}, shows that for $m = 5$, the conjugacy classes of $A_5$ with size $3p^n$ are precisely those stated in the theorem. In contrast, $A_6$, $A_7$ and $A_8$ possess no conjugacy classes of size $3p^n$ for any prime $p$ and any integer $n \ge 0$.

For $m \ge 9$, we use the list of subgroups of $A_m$ whose indices are divisible by exactly two distinct primes, given in  \cite[Table 3.1]{LL}.   According to the lower section of this table, the only  subgroups of $A_m$ with $m \ge 9$ whose indices are powers of two primes  are $A_{m-1}$, $S_{m-2}$, and $A_{m-2}$, subject to certain arithmetic conditions on $m$. However, all these subgroups have a trivial center when $m\geq 9$, making it impossible for them to coincide with the centralizer of any non-trivial element of $A_m$. This proves the non-existence of conjugacy classes of size $3p^n$ for $m \ge 9$.
\end{proof}

\begin{theorem} \label{thm1.2}
No sporadic simple group possesses a conjugacy class of size $3p^n$ for some prime $p$ and integer $n\geq 0$.
\end{theorem}

\begin{proof} 
One can check the result by  using the Atlas of Finite Groups  \cite{Con} or via the GAP Character Table Library \cite{GAP4} by computing the sizes of the conjugacy classes of the sporadic simple groups (and the Tits group), which can be read off from their well-known character tables. Alternatively, the result can be deduced from  \cite[Theorem 1.1 and Table 3.2]{LL}. The only possible pairs $(G,H)$ where $G$ is a sporadic simple group and $H\leq G$ has index $3p^n$ for some prime $p$ and   integer $n\geq 0$, are $(M_{11}, {\rm PSL}_2(11))$, $(M_{22}, M_{11})$, and $(M_{24},M_{23})$. But none of these subgroups can be the   centralizer of any element of $G$,  since all of them have trivial center.
\end{proof}

In the following result, we deal with simple groups of Lie type that possess a conjugacy class of size $3p^n$. We show there exist only  two  such simple groups, one of them being ${\rm PSL}(2,5)\cong A_5$, which already appeared in Theorem \ref{thm1.1}. Consequently, we omit  further details for this group in the statement below.

\begin{theorem}\label{thm1.3}
Let $G$ be a simple group of Lie type with a conjugacy class of  size $3p^n$ for some prime $p$ and integer $n \ge 0$. Then $G\cong  {\rm PSL}(2,5)$ or $G \cong \operatorname{PSL}(2,7)$. Furthermore, in the latter case 
\begin{itemize}
    \item[(a)]
$G$ has a conjugacy class of size $21 = 3 \cdot 7$, consisting of elements of order $2$ lying in the center of the Sylow $2$-subgroups.
\item[(b)] $G$  has a  conjugacy class of size $24 = 3 \cdot 2^3$, consisting of elements of order $7$.
\end{itemize}
\end{theorem}

\begin{proof} Let $x\in G$ with $|x^G|=3p^n$.   Since $G$ is a simple group, then  Kazarin's result \cite[Lemma]{Kaz} allows us to assume that $p\ne 3$. Consequently, the pair $(G,C_G(x))$ must belong to the classification provided in Theorem 1.1 of \cite{LL}. By inspecting \cite[Table 5.1]{LL}, we can rule out all the exceptional groups of Lie type since none of them have subgroups of index $3p^n$. Therefore, the candidates for $(G,C_G(x))$
are listed in  Tables 4.1, 4.2, 4.3 and 4.5 of \cite{LL}. We discuss the pairs appearing according to the family of $G$. 

\smallskip
(1) $G$ is a  linear group. Then the candidates for the pair $(G,C_G(x))$  can only  appear in the lower section of \cite[Table 4.1]{LL} and  comprise the following.

\smallskip

(1.1) $(G,C_G(x))=({\rm PSL}(2,2^{2^i}),D_{2(2^{2^i}-1)})$ or $(G,C_G(x))=({\rm PSL}(2,2^{2^i}),Z_{ 2^{2^i}-1})$. For the first case,  the index is given by $|G:C_G(x)|=2^{2^i-1}(2^{2^i}+1)=3p^n$, which requires that  $(2^{2^i}+1)=3$. It follows that $i=0$ and thus $G={\rm PSL}(2,2)\cong S_3$, contradicting the simplicity of $G$. In the second case, we have $|G:C_G(x)|=2^{2^i}(2^{2^i}+1)=3p^n$, which gives 
 $2^{2^i}+1=3$, leading to the same contradiction.

\smallskip
(1.2) $(G,C_G(x))=({\rm PSL}(2,2^f),D_{2(2^f+1)})$ or $({\rm PSL}(2,2^f),Z_{(2^f+1)})$. In the former case,  the index is given by  $|G:C_G(x)|=2^{f-1}(2^{f}-1)=3p^n$. Then $(2^{f}-1)=3$, so $f=2$, and thus $G={\rm PSL}(2,2^2)\cong A_5$. This case is already dealt with in Theorem \ref{thm1.1}, and gives rise to the first group in the statement of the theorem.  In the latter case, we have $|G:C_G(x)|=2^f(2^f-1)=3p^n$, so we obtain $2^f-1=3$. It follows that $f=2$, and again $G={\rm PSL}(2,4)\cong A_5$.

\smallskip
(1.3) $(G,C_G(x))=({\rm PSL}(2,s^{2i}),D_{s^{2i}-1})$ with $\frac{1}{2}(s^{2i}+1)$ a prime power. Then the index is given by $|G:C_G(x)|=\frac{1}{2} s^{2i}(s^{2i}+1)=3p^n$. These conditions necessarily imply that  $\frac{1}{2}(s^{2i}+1)=3$, and then $s=5$ and $i=0$, giving  again $G={\rm PSL}(2,5)\cong A_5$.

\smallskip
(1.4) $(G,C_G(x))=({\rm PSL}(2,s),D_{s+1})$. Then $|G:C_G(x)|= \frac{1}{2}s(s-1)=3p^n$. The fact that $s$ is prime implies that either  $s=3$, or $s=p$ and $\frac{1}{2}(s-1)=3$. It follows that $s=3$ or $s=7$. As $G$ is a simple group,  $G={\rm PSL}(2,7)$.
In this case, $C_G(x)\cong D_8$, which corresponds to the centralizer of an element $x$ of order 2 lying in the center of a Sylow $2$-subgroup of $G$ \cite{Con}. This yields $\vert{}x^G\vert{}=21=3\cdot 7$, which gives rise to case (a) of the theorem.

\smallskip
(1.5) $(G,C_G(x))=({\rm PSL}(2,2^{2^{i+1}}),{\rm PSL}(2,2^{2^i}))$. In this case, we have $|G:C_G(x)|=2^{2^i}(2^{2^{i+1}}+1)=3p^n$. So we get  $2^{2^{i+1}}+1=3$, which has no integer solution for any $i\geq 0$.

\medskip
 For the remaining cases in \cite[Table 4.1]{LL}, we adopt the standard notation of \cite{KL}. Let $S$ be a Sylow $s$-subgroup of  a simple group of Lie type $G$ defined over a field of characteristic $s$. For a positive integer $i$, let $L_i$ be the Levi subgroup.  Then $B=N_G(S)=SD$ is the Borel subgroup of $G$, and $P_i=BL_i=(O_s(P_i)L_i)D$ are parabolic subgroups with $O_s(P_i)L_i\lhd P_i$, where $D$ is the diagonal subgroup of $G$. The subgroup $L_iD$ is referred to as the Levi complement (or Levi factor) of $P_i$.

\smallskip
(1.6) $(G,C_G(x))=({\rm PSL}(m,q), P_1)$, where $q=s^f$ and $s$ is a prime. In this case, $P_1=C_G(x)$ and has index $(q^m-1)/(q-1)= q^{m-1}+\ldots + q +1$, which is clearly coprime to $s$. Thus, $P_1$ contains a Sylow $s$-subgroup $S$ of $G$, and  we obtain  $x\in O_s(P_1)$ by \cite[Lemma 3.3]{BGL} and $x\in Z(S)$. Following the established notation, we write $P_1=BL_1 = (O_s(P_1)L_1)D$. It is known that $L_1D$ acts irreducibly on $O_s(P_1)$, which, combined  with the fact that $x\in Z(S)$ forces $O_s(P_1)=\langle x\rangle$. Furthermore, from this irreducibility we deduce that $o(x)=s$. Now,  by \cite[Theorem 4.1.17]{KL}, we know  that $|O_s(P_1)|=q^{m-1}=s^{f(m-1)}=s$. Therefore, we conclude that $f=1$ and $m=2$, that is, $G\cong {\rm PSL}(2,s)$. However, since $s\geq 5$, it follows by \cite[Hauptsatz II.8.27]{Hup67} that $P_1=O_s(P_1)\rtimes C_{\frac{s-1}{2}}$ is a Frobenius group, which contradicts the fact that $x\in Z(P_1)$.

\smallskip
 (1.7) $(G,C_G(x))=({\rm PSL}(m,q), H)$, where $q=s^f$ with $s$ prime,   $H$ satisfies $O_s(P_1)H=P_1$, and we have $|G:C_G(x)|=s^a\frac{q^m-1}{q-1}$ for some $a\geq 1$. By our assumption, $s^a\frac{q^m-1}{q-1}=3p^n$, and since the factors on the left-hand side are coprime, we get either   $s^a=3$ or $\frac{q^m-1}{q-1}=3$. In the latter case, we have $3q-q^m= 2$, so $q=2$ and $m=2$, yielding $G\cong {\rm PSL}(2,2)\cong S_3$, a contradiction. As a result, we can assume $s^a=3$, that is, $s=3$ and $a=1$, and hence $\frac{q^m-1}{q-1}=p^n$. We now distinguish between the cases $m=2$ or $m\geq 3$.
 
 If $m=2$, then $q+1=p^n$. Since $q=3^f$, it follows that $p=2$, and hence $3^f=2^n-1$. This forces $n$ to be even.  The case $n=2$ leads to $G\cong {\rm PSL}(2,3)\cong A_4$, which is not simple. Thus, we can write $n=2n_0$ with $n_0>1$, and then $2^n-1=(2^{n_0}+1)(2^{n_0}-1)$.   As $(2^{n_0}+1,2^{n_0}-1)=1$, we get a contradiction. 
 
 We consider the case $m\geq3$. Note that $H\leq P_1$ and $O_s(P_1)\unlhd P_1$, which implies that $H\cap O_s(P_1)\unlhd H$. Since $O_s(P_1)\cong ({\mathbb F}_q^{m-1},+)$ is abelian, it follows that  $H\cap O_s(P_1)\unlhd HO_s(P_1)=P_1$. On the other hand,  $P_1=O_s(P_1)L_1D$ where $L_1D$ acts irreducibly on $O_s(P_1)$. Then we obtain $H\cap O_s(P_1)=1$. Consequently, $$|O_s(P_1)|=|O_s(P_1)H: H |= |P_1: H|=s^a=3.$$
 Now, by \cite[Theorem 4.1.17]{KL}, we know that $|O_s(P_1)|=q^{(m-1)}=3^{f(m-1)}$.  We conclude that $f=1$ and $m=2$, which gives $G={\rm PSL}(2,3)\cong A_4$,  a contradiction.

\smallskip
  (1.8) Suppose that $(G,C_G(x))=({\rm PSL}(m,q), O_s(P_1)L_1R)$, where $q=s^f$ with $s$ prime, and $R<D$ with $|D:R|=r^b\mid q-1$.
In addition, we have $|G:C_G(x)|=r^b\frac{q^m-1}{q-1}=3p^n$. It is straightforward to see that $s\nmid r^b\frac{q^m-1}{q-1}$, from which we deduce that $x$ centralizes a Sylow $s$-subgroup of $G$. By \cite[Lemma 3.3]{BGL}, it follows that $x\in  O_s(P_1)$. Arguing as in case (1.6), we obtain $G\cong {\rm PSL}(2,s)$. If $L_1R>1$, then $O_s(P_1)L_1R$ is a Frobenius group. On the other hand, we have $C_G(x)=O_s(P_1)L_1R$, leading to a contradiction. This allows us to assume that $L_1=R=1$, so that $C_G(x)=O_s(P_1)$. Then the equality $|G:C_G(x)|=3p^n$ implies that $|G|=3\cdot s\cdot p^n$ (with $r=3$). By  applying \cite[Theorem 4]{Her}, which classifies finite simple groups with exactly three prime divisors, we obtain $G\cong {\rm PSL}(2,5), {\rm PSL}(2,7)$ or ${\rm PSL}(2,17)$. However, the parabolic structure requires $r\mid s-1$ with $r=3$, so this forces $G\cong {\rm PSL}(2,7)$. In this case, an element $x\in G$ of order $7$ yields $|x^G|=|G:C_G(x)|=24=3\cdot 2^3$, giving case (b) of the statement.

\smallskip
(1.9) If $(G,C_G(x))=({\rm PSL}(4,2^{2^i}), P_1)$, then $|G:C_G(x)|=(2^{2^i}+1)(2^{2^{i+1}}+1)=3p^n$. This requires that $2^{2^i}+1=3$ because $2^{2^{i+1}}+1\geq5$, so $i=0$ and $G\cong {\rm PSL}(4,2)\cong A_8$. In this case (which was also covered in Theorem \ref{thm1.1}), we know that  $P_1=O_2(P_1)L_1\cong 2^3:{\rm PSL}(3,2)$, a subgroup with trivial center, which yields a contradiction.

\smallskip
(2) Suppose that  $G$ is a symplectic group. According to \cite[Table 4.2]{LL}. it suffices to inspect its lower section. In almost all cases, the index of the subgroup cannot be of the form $3p^n$, being the only exceptions the pairs $({\rm Sp}_4(s),P_1)$ and $({\rm Sp}_{2m}(2^f), {\rm O}_{2m}^-(2))$.

\smallskip
(2.1) Assume first $(G,C_{G}(x))=(\mathrm{Sp}_{4}(s),P_{1})$. In this case, we have $|G:C_{G}(x)|=(s+1)(s^{2}+1)=3p^{n}$. 
Note that $(s+1, s^{2}+1) = (s+1,(s-1)(s+1)+ 2)=(s+1,2) \le 2$. 
If $s$ is odd, then $(s+1, s^{2}+1)=2$, forcing $4$ to divide $|g^G|$, and hence $p=2$. Then $(s+1)(s^{2}+1)=3\cdot 2^{n}$, which forces  one of the factors, either $s+1$ or $s^2+1$, to be a power of $2$. However, neither $s^2+1=2^a$ for some $a\geq 0$, nor $s+1=2^b$ with $s^2+1=3\cdot 2^c$ for some $b,c\geq 1$ have integer solutions given that every odd prime $s$ satisfies  $s^2 +1\equiv 2 \pmod 4$.
Therefore, we must have $s=2$, which gives $G= \mathrm{Sp}_{4}(2)\cong S_{6}$, contradicting the simplicity of $G$.

\smallskip
(2.2) Next consider $(G,C_G(x))=({\rm Sp}_{2m}(2^f), {\rm O}_{2m}^-(2))$, where $m$ is a prime and $f=1$. We have $|G:C_G(x)|=2^{m-1}(2^m-1)=3p^n$, which implies that $2^m-1=3$, so $m=2$.  As a consequence, $G= {\rm Sp}_4(2)\cong S_6$,  the same contradiction as in case (2.1).

\smallskip
(3)  Assume that $G$ is a  unitary group. Looking at \cite[Table 4.3]{LL}, it suffices to consider the pairs $(G,H)$ appearing in the lower section of the table. We show that all these cases can be ruled out. First, note that $(s+1,s^2-s+1)$ divides $3$. Indeed, $(s+1,s^2-s+1)=(s+1, (s+1)(s-2)+3)=(s+1,3)$, which clearly divides $3$. On the other hand, according to \cite[Table 4.3]{LL}, we observe that both factors, $(s+1)$ and $(s^2-s+1)$,  divide the index $|G: C_G(x)|= 3p^n$, with $p\neq3$ and $n \geq 1$. This forces $(s+1,s^2-s+1)=1$, which in particular gives $s\neq2$. However, since one of the factors $s^2-s+1$ and $s+1$ must be equal to $3$,  both possibilities trivially lead to a  contradiction.

\smallskip
(4) Finally, suppose that $G$ is an  orthogonal group. Following  the cases listed in \cite[Table 4.4]{LL}, it follows that the only candidate for the pair $(G,C_G(x))$ is $({\rm \Omega}_{2m+1}(q),Z_2\times \Omega_{2m}^+(q))$, where $q=s^f$ for a prime $s$. In addition, in this case $|G:C_G(x)|=\frac{1}{2}q^m(q^m+1)=3p^n$, where $q^m=2r^a-1$ and $r$ is a prime. In particular, $q$ is odd, which forces either  $q^m=3$ or $\frac{1}{2}(q^m+1)=3$. This gives either $q=3$ and $m=1$, or $q=5$ and $m=1$. Consequently, $G\cong \Omega_3(3)$ or $\Omega_3(5)$. However, $\Omega_3(3)\cong {\rm PSL}(2,3)$ is not simple, while $\Omega_3(5)\cong {\rm PSL}(2,5)$ has already been considered. The remaining orthogonal group, $\Omega_{2m+1}(3)$, listed in \cite[Table 4.5]{LL}, is also ruled out because the arithmetic conditions force $m=1$.
\end{proof}

\section{Proof of Theorem B}

The proof of Theorem B relies on analysing the indices of fixed-point subgroups of automorphisms acting on non-abelian finite simple groups. 
Recall that Kazarin \cite[Corollary 1]{Kaz} proved, without employing Classification,  that a non-trivial automorphism of a non-abelian finite simple group cannot fix a subgroup of prime-power index. To extend this study to indices of the form $3p^n$, our analysis is naturally divided into three main cases according to the CFSG: 
 alternating groups (Theorem \ref{thm3.1}), sporadic groups (Theorem \ref{thm3.2}), and groups of Lie type (Theorem \ref{thm3.3}). Once these  results are established, we combine them at the end of the section to complete the proof of Theorem B.

\begin{theorem}\label{thm3.1}
Let $G$ be the alternating group $A_m$ with $m\geq 5$ and $\alpha \in \operatorname{Aut}(G)$ such that 
$|G : C_G(\alpha)| = 3p^n$
for some prime $p$ and integer $n \ge 0$. Then one of the following holds:
\begin{itemize}
\item[(a)] $G=A_5$ and  $\alpha$ is an inner automorphism  of order $2$ or $5$  and  $|G : C_G(\alpha)| = 15$ or $12$, respectively. 

\item[(b)] $G=A_6$ and $\alpha$ is  an  automorphism induced by either a transposition of $S_6$ or the product of three disjoint transpositions of $S_6$. In both cases, $|G : C_G(\alpha)| = 15$.
\item[(c)] $G=A_7$ and $\alpha$ is an automorphism induced by a transposition of $S_7$. In this case, $|G : C_G(\alpha)| = 21$.
\end{itemize}
\end{theorem}

\begin{proof}

   We examine the index $|G : C_G(\alpha)| = 3p^n$ for distinct values $m$.    
     In view of Theorem \ref{thm1.1}, the group $G = A_5$ satisfies the condition when $\alpha$ is an inner automorphism of order $2$ or $5$. In the former case, $p=5$ and $n=1$, and the index is $15 = 3 \cdot 5$; in the latter case $p=2$ and $n=2$, and the index is $12 = 3 \cdot 2^2$. On the other hand, since ${\rm Aut}(G)\cong S_5$, we can easily compute the fixed point subgroup in $A_5$ of any odd permutation in $S_5$, and easily follows that its index in $A_5$ is not of the form $3p^n$. Therefore,
     we obtain (a).

 For $G = A_6$, recall that ${\rm Out}(A_6)\cong C_2\times C_2$. By Theorem \ref{thm1.1}, we know that $\alpha$ cannot be an inner automorphism, so we assume  $\alpha\in {\rm Aut}(A_6)\setminus {\rm Inn}(A_6)$. We know that $\Gamma=G\langle \alpha\rangle $ is isomorphic to one of the three possible extensions of degree $2$ of $A_6$ within ${\rm Aut}(A_6)$, namely, $S_6$, ${\rm PGL}(2,9)$ and $M_{10}$. Note that $|\Gamma:C_{\Gamma}(\alpha)|=|G: C_G(\alpha) |$, so it is enough to compute the class sizes of each one of these extensions (with \cite{Con} or \cite{GAP4}). We find that the only possibility is $\Gamma=S_6$ and  $|A_6 : C_{A_6}(\alpha)| =  15 = 3 \cdot 5$, so $p = 5$ and $n = 1$. Furthermore, this occurs if and only if $\alpha$ belongs to the conjugacy class of transpositions or to the class formed by the product of three disjoint transpositions in $S_6$. This gives case (b).
 
  For $G = A_7$,  
 we have ${\rm Out}(A_7)\cong C_2$, and the only nontrivial extension of $A_7$ within ${\rm Aut}(A_7)$ is $S_7$. Computing the class sizes of $S_7$, we have $|S_7: C_{S_7}(\alpha)|= 21=3\cdot 7 $, when $\alpha $ is any transposition of $S_7$. Also, we have $C_{A_7}(\alpha)\cong S_5$ and $|A_7 : C_{A_7}(\alpha)| = 21$. In fact, transpositions induce the only automorphisms that satisfy the required condition. Thus, we obtain (c).

Suppose $G=A_8$. Similarly to  $A_7$,  the reasoning reduces  to computing the class sizes of $S_8$, and we find that there do not exist conjugacy  classes of size $3p^n$.

 Finally, for $G=A_m$ with  $m \ge 9$ we also have ${\rm Aut}(G)=S_m$. We prove that no automorphism  satisfies the class size condition in the statement. First, we show that  $\alpha$ can be  assumed to be of prime order.
 Suppose that $o(\alpha) = k > 1$ is not prime, and let $r$ be a prime factor of $k$. Consider  $1\neq \beta=\alpha^{k/r} \in {\rm Aut}(G)$, which has prime order $r$. We have $C_G(\alpha) \le  C_G(\beta) \le G$. Consequently, $|G :C_G(\beta)|$ divides $|G : C_G(\alpha)| = 3p^n$, so $|G :C_G(\beta)| = p^m$  or$|G :C_G(\beta)| = 3p^m$ for some $0 \le m \le n$. 
However, it is impossible that $|G : C_G(\beta)| = p^m$ or $3$, since Corollary 1 of \cite{Kaz} asserts that a non-trivial automorphism  cannot fix a subgroup of prime-power index in a non-abelian simple group.  Therefore, $|G :C_G(\beta)| = 3p^m$ with $m\geq 1$, so we can replace $\alpha$ by $\beta$.
 
 Also, by Theorem \ref{thm1.1}, $\alpha$ cannot be an inner automorphism, and this means that $\alpha$ is an odd permutation in $S_m$. Notice that the only elements of prime order in $S_m\setminus A_ m$  are the products of an odd number $t$ of disjoint transpositions with $2t\leq m$. The size of such conjugacy class in $S_m$ are given by the standard formula $$|\alpha^{S_m}|= \frac{m!}{t!2^t (m-2t)!} =\binom{m}{2t} \frac{(2t)!}{t!2^t}=\binom{m}{2t}(2t-1)!!,$$
 where $(2t-1)!!$ denotes the product of all odd numbers from $1$ up to $2t-1$.  Taking into account that $m\geq 9$,  next we check that $|\alpha^G|$ cannot be equal to $3p^n$. Indeed, when $t=1$, the equation reduces to $m(m-1)=6p^n$, but since $m$ and $m-1$ are coprime integers, the unique solutions are $m=4$,  $m=6$, and $m=7$, against our assumption $m \geq 9$. If $t=2$, the equation becomes $\binom{m}{4}=p^n$, which is clearly impossible because $m\geq 9$.  For $t=3$, we obtain the equation $\binom{m}{6}15=3p^n$, leading to  $\binom{m}{6}=5^{n-1}$, again a contradiction. Finally, if $t\geq 4$, then 
$(2t-1)!!$ is certainly divisible by at least three different primes, and hence $|\alpha^{S_n}|$ has at least three distinct prime divisors as well. This proves that there are no automorphisms satisfying the desired condition, as required. 
 \end{proof}

\begin{theorem} \label{thm3.2}
     No sporadic simple group $G$  possesses an automorphism  $\alpha\in {\rm Aut}(G)$ such that $|G:C_G(\alpha)|= 3p^n$, for some prime $p$ and $n\geq 0$.
     \end{theorem}

\begin{proof}
 
For the 26 sporadic simple groups (and the Tits group), an explicit check in \cite[Table 3.2]{LL} shows that the only groups having subgroups of index $3p^n$ for some prime $p$ are the three Mathieu groups $M_{11}$, $M_{12}$ and $M_{24}$. Since these groups have trivial outer automorphism groups, then all automorphisms are inner and then Theorem \ref{thm1.2} shows that there are no sporadic groups possessing an automorphism $\alpha$ such that $|G : C_G(\alpha)| = 3p^n$. 
\end{proof}

\begin{theorem}\label{thm3.3} Let $G$ be a non-abelian finite simple group of Lie type and $\alpha\in {\rm Aut}(G)$ such that $|G:C_G(\alpha)|=3p^n$ for some prime $p$ and integer $n\geq 0$. Then one of the following holds:

\begin{itemize}
\smallskip
\item[(a)]  $G= {\rm PSL}(2,5)$ and   $\alpha$ is an inner automorphism of order $2$ or $5$ with $|G:C_G(\alpha)|=15$ or $12$, respectively;

\smallskip
\item[(b)]  $G= {\rm PSL}(2,7)$ and $\alpha$ is an inner automorphism of order $2$ (corresponding to an element lying in the center of Sylow $2$-subgroup of $G$), or order $7$ with $|G:C_G(\alpha)|=21$ or $24$, respectively;

\smallskip
\item[(c)] $G={\rm PSL}(2,9)$ and $\alpha$ is a non-inner automorphism of order $2$ with  $|G:C_G(\alpha)|=15$.
\end{itemize}
\end{theorem}

 \begin{proof} Assume first that $\alpha\in {\rm Inn}(G)$.  Then the action of $\alpha$  by conjugation corresponds to that of an element $x\in G$ such that $C_G(\alpha) = C_G(x)$. The conjugacy classes of size $3p^n$ in simple groups of Lie type have been classified in Theorem 2.3. This immediately yields cases (a) and (b) of the statement.  Therefore, we will assume for the rest of the proof that $\alpha$ is non-inner.

We now analyze the three groups appearing in the statement.  The group ${\rm PSL}(2,5)\cong A_5$ has been dealt with in Theorem \ref{thm3.1}, and we know that it does not have any non-inner automorphism satisfying the hypothesis. Suppose now that $G={\rm PSL}(2,7)$. We then know that ${\rm  Out}(G)\cong C_2$. Assume further that  there exists $\alpha\in {\rm Aut}(G)\setminus {\rm Inn}(G)$ such that $|G:C_G(\alpha)|=3p^n$. Then $G\langle\alpha \rangle\cong {\rm PGL}(2,7)$ and  it is enough to compute the class sizes of  ${\rm PGL}(2,7)$ (\cite{Con} or \cite{GAP4}) to find that no element out of ${\rm PSL}(2,7)$ has a conjugacy class of size $3p^n$.  Finally, since ${\rm PSL}(2,9)$ is isomorphic to $A_6$, this case has been analyzed in Theorem \ref{thm3.1}, leading directly to case (c). Thus, for the rest of the proof, we  exclude these three simple groups.

 \smallskip
The  strategy of the proof consists of verifying that none of the remaining   simple groups,  which appear in Tables 4.1, 4.2, 4.3, 4.4 and 5.1 of \cite{LL}, gives rise to a non-inner automorphism satisfying the hypothesis of the statement.
For that purpose, we do an inspection of the proof of Theorem \ref{thm1.3}, and this shows that for all  exceptional simple groups of Lie type  as well as the remaining classical groups (symplectic, unitary, and orthogonal groups, along with cases (1.1)–(1.3), (1.5), (1.7), and (1.9) in the case of linear groups), no specific structural information about $C_G(\alpha)$ is needed to discard these groups beyond the arithmetic fact that  $[G:C_G(x)]=3p^n$. Therefore, all these cases can be excluded, and it only remains to consider the following pairs $(G,C_G(\alpha))$ within the family of linear groups: 

\smallskip

(1) $({\rm PSL}(2,19), A_5)$, which appears in the upper section of \cite[Table 4.1]{LL}; 

\smallskip
\noindent
and the following three cases given  in the lower section of \cite[Table 4.1]{LL} that correspond to cases (1.4), (1,6) and (1.8) in the proof of Theorem \ref{thm1.3}:

\smallskip
(2) $({\rm PSL}(2,s)), D_{s+1})$, where $s$ is a prime.

\smallskip
(3) $({\rm PSL}(m, q), P_1)$, with $q = s^f$ and $s$ a prime.  
\smallskip

(4) $({\rm PSL}(m, q),  O_s(P_1)L_1R)$, where $r^b \frac{q^m-1}{q-1} = 3p^n$, $q = s^f$ with $s$ a prime and $r\mid q-1$.  

\smallskip
Next, we discuss them case by case. 

\smallskip
(1) If $(G,C_G(\alpha))=({\rm PSL}(2,19), A_5)$, then  ${\rm PGL}(2,19)$ would be written as  ${\rm PGL}(2,19)={\rm PSL}(2,19)\langle \alpha \rangle$, with $C_G(\alpha)\cong A_5$.  However, by inspecting  \cite{Con} (or using \cite{GAP4}), we find that $A_5$ is not an $\langle \alpha \rangle$-invariant subgroup, yielding a contradiction. 

\smallskip
(2)  $(G,C_G(\alpha))=({\rm PSL}(2,s), D_{s+1})$, where $s$ is a prime.
In this case, ${\rm PGL}(2,s)={\rm PSL}(2,s)\langle \alpha\rangle$. Following the same argument as in case (1.4) of the proof of Theorem \ref{thm1.3}, we obtain that $s=7$, so $G={\rm PSL}(2,7)$ and $C_G(\alpha)\cong D_8$. Let $D=C_{{\rm PGL}(2,7)(}\alpha)\cong D_8 \langle \alpha\rangle$. Notice that $Z(D_8)\leq Z(D)$ and $\alpha\in Z(D)$, so $|Z(D)| \geq 4$. However,  consulting \cite{Con} (or \cite{GAP4}), we have $C_{{\rm PGL}(2,7)}(\alpha) =D_{16}$. This is clearly a contradiction. 

\smallskip
(3) and (4). To deal with these cases,  we first need to do the following reduction. We show that, without loss of generality, we may assume that $\alpha$ has prime order. Suppose that $o(\alpha) = k > 1$ is not a prime number, and let $r$ be a prime factor of $k$. Consider $\beta = \alpha^{k/r} \in {\rm Aut}(G)$, which has prime order $r$. We have $C_G(\alpha) \le  C_G(\beta) \le G$. Consequently, $|G :C_G(\beta)|$ divides $|G : C_G(\alpha)| = 3p^n$, so $|G :C_G(\beta)| = p^m$  or $|G :C_G(\beta)| = 3p^m$ for some $0 \le m \le n$. 
But \cite[Corollary 1]{Kaz} implies that  $|G : C_G(\beta)|$ cannot be equal to $p^m$ or $3$. Thus, $|G :C_G(\beta)| = 3p^m$ with $m\geq 1$. Now, if $\beta$ were an inner automorphism of $G$, then Theorem \ref{thm1.3} would imply that $G\cong {\rm PSL}(2,5)$ or ${\rm PSL}(2,7)$. However, both groups have already been excluded at the beginning of the proof. This proves that $\beta $ can be assumed to be non-inner, and henceforth, by replacing $\alpha$ by $\beta$,  we can assume that $\alpha$ has prime order, as desired.

\smallskip 
   Let $T={\rm PGL}(m, q)$. Then $T\geq {\rm PSL}(m, q)=G$.  Since $q=s^f$, we have $\frac{q^m-1}{q-1}=1+q+\cdots+q^{m-1}\equiv 1\pmod{s}$. Moreover, the fact that $r\mid(q-1)$ implies that $(r,s)=1$. Thus, as $|G:C_G(\alpha)|=\frac{q^m-1}{q-1}$  or $r^b\frac{q^m-1}{q-1}$ in cases (3) and (4), respectively, the index $|G:C_G(\alpha)|$ is coprime to $s$. Consequently, $C_G(\alpha)$ contains a Sylow $s$-subgroup of $G$.  This condition plays a crucial role to discard these cases. Indeed, since $\alpha$ is a prime-order  automorphism of $G$, we can determine the order of $C_T(\alpha)$ using \cite[Proposition 3.2.16 and Table B.3]{BG}.  Specifically, it suffices to compare the $s$-part $|C_T(\alpha)|_s$ with the $s$-part $|G|_s=q^{m(m-1)/2}$. We will demonstrate that $|C_T(\alpha)|_s < |G|_s$, so this will contradict that $C_G(\alpha)$ contains a Sylow $s$-subgroup of $G$. 

\smallskip
We distinguish the various types of  prime-order automorphisms following the classification in \cite[Table B.3]{BG};  namely, $\texttt{u}$ (inner), $\texttt{s}$ (diagonal),  $\texttt{f}$  (field), $\texttt{gf}$ (graph-field) and $\texttt{g}$ (graph).  Type ${\texttt u}$ is excluded because $\alpha$ is non-inner, so we proceed with a case-by-case analysis for the remaining types:

\smallskip
(3.1)  $\alpha$  is of type $\texttt{s}$ with $i > 1$.  Note that 
$$|C_T(\alpha)| = \frac{1}{q-1} |\mathrm{GL}(e, q)| \prod_j |\mathrm{GL}(a_j, q^i)|$$
(where $e$ and $a_j$ are defined as in \cite[Table B.3]{BG}). Using standard formulas for the orders of general linear groups, we have
$$|\mathrm{GL}(e, q)| = q^{\frac{e(e-1)}{2}} \prod_{t=1}^e (q^t - 1)\,\, {\rm and} \, \,|\mathrm{GL}(a_j, q^i)| = q^{\frac{a_j(a_j-1)i}{2}} \prod_{t=1}^{a_j} (q^{it} - 1).$$
Now, since $m = e + i \sum_j a_j$, the $s$-part of $|C_T(\alpha)|$ evaluates to
 $|C_T(\alpha)|_s = q^{\frac{e(e-1) + i \sum_j a_j(a_j-1)}{2}}$.
Given that $i > 1$, we have
 $$e(e - 1) + i \sum_j a_j(a_j - 1) < \left(e + i \sum_j a_j\right)^2 - \left(e + i \sum_j a_j\right) = m^2 - m,$$
which guarantees that $|C_T(\alpha)|_s < q^{\frac{m(m-1)}{2}} = |G|_s$. 

\smallskip
(3.2)  $\alpha$  is of type $\texttt{s}$ with $i=1$. 
Then 
$$| C_T(\alpha)|=\frac{1}{q-1} \prod_{j}|{\rm GL}(a_j,q)|,$$ with $\sum_{i} a_i=m$.
Since $|C_T(\alpha)|_s= \prod_{j}|{\rm GL}(a_j,q)|_s=q^{\sum_j a_j(a_j-1)/2}$ and $\sum_{j} a_j=m$, it follows that $\sum_j a_j(a_j-1)/2<m(m-1)/2$.  We conclude  that  $|C_T(\alpha)|_s< |G|_s$.  

\smallskip
(3.3)  $\alpha$  is of type $\texttt{s}$ with $i=1$, the order of $\alpha$ divides $m$ and  $| C_T(\alpha)|=\frac{r}{q-1} |{\rm GL}(\frac{m}{r},q)|^r$.   
Since  $|{\rm GL}(\frac{m}{r},q)|_s^r=q^{\frac{m}{r}(\frac{m}{r}-1)r/2}$,  we get   $| C_T(\alpha)|_s< |G|_s$.

\smallskip
(3.4)  $\alpha$  is of type $\texttt{s}$ with $i=1$, the order of $\alpha$  divides $m$ and  $| C_T(\alpha)|=\frac{r}{q-1} |{\rm GL}(\frac{m}{r},q^r)|$.  
Proceeding as in case (3.3) we obtain the same conclusion.

\smallskip
(3.5)  $\alpha$  is of type $\texttt{f}$ with  $o(\alpha)\mid f$. 
Again this case is similar to (3.3). 

\smallskip
(3.6)  $\alpha$  is of type $\texttt{gf}$ with $o(\alpha)=2$, $f$ even and $m\geq3$.
Then  $|C_T(\alpha)|= |{\rm PGU}(m,q^{1/2})|$.  Since  $$|{\rm PGU}(m,q^{1/2})|=q^{\frac{1}{2} (\frac{m(m-1)}{2})}\frac{1}{q+1}\prod_i^m(q^i-(-1)^i),$$  we obtain $|C_T(\alpha)|_s=q^{\frac{m(m-1)}{4}}$.  Hence  $|C_T(\alpha)|_s< |G|_s$.

\smallskip
(3.7)  $\alpha$  is of type $\texttt{g}$ with $o(\alpha)=2$  and $m$   odd.
This case can be  handled   similarly to (3.3).

\smallskip
(3.8)  $\alpha$  is of type $\texttt{g}$ with $o(\alpha)=2$ and $m\geq4$ even.
Likewise, this case follows by calculation similar to that in (3.3).

\smallskip
(3.9)  $\alpha$  is of type $\texttt{g}$ with $o(\alpha)=2$,  $m\geq4$ even and $s\neq2$.
In addition, $C_T(\alpha)={\rm PGO^+}(m,q)$. Since  
$$|{\rm PGO^+}(m,q)|=q^{m(m/2-1)/2}(q^{m/2}-1)\prod_{i=1}^{\frac{m}{2}-1}(q^{2i}-1).$$ we get $|C_T(\alpha)|_s=q^{m(m-2) /4}$. Hence $|C_T(\alpha)|_s< |G|_s$.

\smallskip
(3.10)  $\alpha$  is of type $\texttt{g}$ with $o(\alpha)=2$,  $m\geq4$ even and $s\neq2$. In addition, $C_T(\alpha)={\rm PGO^-}(m,q)$. This case is analogous to (3.9).

\smallskip
In conclusion, in all cases we have shown that  $|C_T(\alpha)|_s<|G|_s$. Since $|C_G(\alpha)|_s\leq |C_T(\alpha)|_s$, it follows that $|C_G(\alpha)|_s<|G|_s$, which provides the required contradiction.        
 \end{proof}

For proving Theorem B, we also make use of following result of elementary character, which we would like to state, in case it is not widely known.

\begin{lemma} \label{cent}
Suppose that a finite group $A$ acts on a finite group $G$ which allows
a direct decomposition $G = H_1\times \cdots \times H_n$, that is invariant under $A$, i.e.,
$H^a_i\in \{H_1,...,H_n\}$ for all $a \in A$ and all $i \in \{1,...,n\}$. Assume further that
A acts transitively on $\{H_1,...,H_n\}$. Let $H \in  \{H_1,...,H_n\}$ and $B = {\bf N}_A(H)$. Then ${\bf C}_G(A)\cong {\bf C}_H(B)$.
\end{lemma}

\begin{proof}
This is \cite[8.1.6]{KS}.
\end{proof}

We are ready to prove Theorem B, which we state again.

\begin{theorem}
Let $G$ be a finite group having a conjugacy class $g^G$ of size $3p^n$, where $p$ is prime and $n \ge 1$, such that $\langle g^G \rangle$ is non-solvable.  Let $\hat{g}$ denote the image of $g$ in $\langle g^G\rangle/S(\langle g^G\rangle)$. Then  $\langle g^G\rangle/S(\langle g^G\rangle)$ is isomorphic to one of the following  groups:
\begin{enumerate}
    \item $A_5$, where $\hat{g}$ is an element of order $2$ or $5$.
    \item $S_6$, where $\hat{g}$ is an odd permutation of order $2$ in $S_6$ (either a transposition or a product of three disjoint transpositions).
    \item $S_7$, where $\hat{g}$ is a transposition in $S_7$.
    \item $\mathrm{PSL}(2, 7)$, where $\hat{g}$ is an element of order $7$, or an element of order $2$ lying in the center of a Sylow $2$-subgroup of $\mathrm{PSL}(2, 7)$.
   \end{enumerate}
\end{theorem}

\begin{proof}
Let $S(G)$ denote the solvable radical of $G$, and set 
$$\overline{G} = G/S(G) \quad {\rm and} \quad \overline{N} = \langle g^G \rangle S(G)/S(G).$$ Throughout the proof, for any  $x \in G$, we denote by $\overline{x}$ its canonical image in  $\overline{G}$. Noting that $S(G)\cap \langle g^G \rangle =S(\langle g^G\rangle)$, the  isomorphism $\langle g^G \rangle / S(\langle g^G \rangle) \cong \overline{N}$ allows us to identify both quotients, and likewise the element $\hat{g}$, appearing in the statement, with $\overline{g}$.

Since $\langle g^G \rangle$ is non-solvable, it follows that $\overline{N} \neq 1$. Let $\overline{L} = \text{Soc}(\overline{G}) = L_1 \times \cdots \times L_k$ be the socle of $\overline{G}$, where each $L_i$  is a minimal normal subgroup of $\overline{G}$, which, in addition, is a direct product of non-abelian simple groups for $i = 1, \dots, k$.

We first prove that $\overline{N}$ centralizes all factors of $\overline{L}$ except one. It is clear that $C_{\overline{G}}(\overline{L}) = 1$, so $\overline{N}$ cannot centralize $L_i$ for all $i$. Assume, for instance, that $L_1$ is not centralized by $\overline{N}$. This implies that $\overline{g}$ does not centralize $L_1$ either; otherwise, all conjugates of $\overline{g}$ in $\overline{G}$ would centralize $L_1$, thereby $\overline{N}$ as well. On the other  hand, as $\overline{g}$ normalizes every $L_i$, we have $ C_{\overline{L}}(\overline{g}) = C_{L_1}(\overline{g})\times\ldots\times C_{L_k}(\overline{g})$, and thus
$$1 \neq |\overline{L} : C_{\overline{L}}(\overline{g})| = \prod_{i=1}^k |L_i : C_{L_i}(\overline{g})|.$$
Note that this index divides $|\overline{G}:C_{\overline{G}}(\overline{g})|$, which divides $|g^G| = 3p^n$.  As a result, $|L_i : C_{L_i}(\overline{g})|$ divides $3p^n$ for every $i$. But these indices cannot be a prime power by \cite[Corollary 1]{Kaz}, so all of them are equal to $1$ except one, say $|L_1 : C_{L_1}(\overline{g})|$, which must be equal to $3p^a$ for some $1\leq a \leq n$. Therefore, $\overline{g}$ centralizes $L_i$ for all $i = 2, \dots, k$, as required.

\medskip
Write $L_1=S_1\times \ldots \times S_r$, where the $S_i$ are isomorphic non-abelian simple groups. Next we prove that the action of $\overline{g}$ on $\{S_1, \dots, S_r\}$, does not have any orbit of length $t \geq 2$. Suppose, for instance, that $\{S_1, \dots, S_t\}$ is the orbit of $S_1$ under the action of $\overline{g}$, and set $T=S_1\times\ldots\times S_t$. By  applying Lemma \ref{cent} to the action of $\overline{g}$ on $T$, we obtain $C_{T}(\overline{g}) \cong C_{S_1}(N_{\langle \overline{g} \rangle}(S_1)) = C_{S_1}( \overline{g} ^s )$ for some $s\geq 1$. Then
$$1\neq |T : C_T(\overline{g})| =\frac{|S_1\times \ldots \times S_t|}{|C_{S_1}(\overline{g}^s)|}= |S_1 : C_{S_1}(\overline{g}^s)| \prod_{j=2}^t |S_j|.$$
Moreover, since $T$ is subnormal in $\overline{G}$, this index divides $|\overline{G}:C_{\overline{G}}(\overline{g})|$, which in turn divides the class size $|G:C_G(g)|=3p^n$. Now, if $t>1$, the factor $|S_1|^{t-1}$ introduces at least three distinct prime divisors in $|G:C_G(g)|=3p^n$, which is impossible. Therefore, $t = 1$, which means that  $\overline{g}$ normalizes $S_i$ for all $i=1,\ldots,r$. We then have 
$$1\neq |{L_1} : C_{L_1}(\overline{g})| = \prod_j |S_j : C_{S_j}(\overline{g})|.$$ If $\overline{g}$ acted non-trivially on two or more simple factors, say $S_1$ and $S_2$, then the product of the corresponding centralizer indices, which cannot be prime powers as mentioned above,  would have more than two factors containing $3$, leading to the fact that $9$ must divide $|g^G|$, a contradiction. We conclude that $\overline{g}$ acts non-trivially on exactly one simple factor, say $S_1$,  and trivially on $S_j$ for $2\leq j \leq r$.
\medskip

We have shown above that $\overline{N}$ centralizes all $L_i$ for $2\leq i \leq k $ and all  $S_j$ for $2\leq j \leq r$. Consequently, $\overline{N} \cap L_j \le Z(L_j) = 1$ for $2\leq j \leq k$. It follows that $L_1$ is the only factor of $\overline{L}$ that can be contained in $\overline{N}$; in fact, as $L_1$ is minimal normal in $\overline{G}$, we necessarily have $\overline{N}\cap \overline{L}=L_1$. Now, if $r>1$, it would follow that $L_1$ centralizes $S_j$ for all $2\leq j\leq r$, which is impossible. Therefore, $r=1$, that is,  $L_1=S_1$. In particular, $S_1$ is normal in $ \overline{N}$. Moreover, we have 
 $C_{\overline{N}}(S_1) \leq C_{\overline{N}}(\overline{L}) \leq C_{\overline{G}}(\overline{L}) = 1$, and hence $\overline{N}$ embeds into
 $\text{Aut}(S_1)$. Thus,
$$S_1 \leq \overline{N} \leq \text{Aut}(S_1).$$
On the other hand, since $S_1$ is subnormal in $\overline{G}$, then $|S_1 : C_{S_1}(\overline{g})|$ divides $|\overline{G} : C_{\overline{G}}(\overline{g})|$, which in turn divides $|g^G| = 3p^n$.  By  \cite[Corollary 1]{Kaz}, we know that $\vert{}S_1 : C_{S_1}(\overline{g})\vert{}$ cannot be a prime power, and as a consequence,  $|S_1 : C_{S_1}(\overline{g})|=3p^b$ with $1\leq b\leq n$. Applying Theorems 3.1--3.3 to the non-trivial action of $\overline{g}\in {\rm Aut}(S_1)$ on the simple group $S_1$, we obtain one of the following cases:

\begin{enumerate}
    \item $S_1 \cong A_5$, where $\overline{g}$ acts as an inner automorphism of order $2$ or $5$, so $\overline{N} \cong A_5$;
    \item $S_1 \cong A_6$, where $\overline{g}$ corresponds to an odd permutation of order $2$ in $S_6$ (either a transposition or a product of three disjoint transpositions), so $\overline{N} \cong S_6$;
    \item $S_1 \cong A_7$, where $\overline{g}$ corresponds to a transposition in $S_7$, so $\overline{N} \cong S_7$;
    \item $S_1 \cong \mathrm{PSL}(2,7)$, where $\overline{g}$ is an inner automorphism of order $7$ or an element of order $2$ lying in the center of a Sylow $2$-subgroup of $\mathrm{PSL}(2,7)$, so $\overline{N} \cong \mathrm{PSL}(2,7)$;
    \end{enumerate}

Finally, since $\overline{N} \cong \langle g^G \rangle / S(\langle g^G \rangle)$, identifying $\overline{g}$ with $\hat{g}$ yields each of the corresponding cases in the statement of the theorem. The proof is complete.  
\end{proof}

\begin{corollary}
 Let $G$ be a finite group having a conjugacy class $g^G$ of size $3p^n$ for some prime $p$ and $n\geq 1$. If $(o(g), 70)=1$, then $\langle g^G \rangle$ is solvable.
 
\end{corollary}

\begin{proof} 
Since the order of $\overline{g}$ in $\langle g^{G}\rangle/S(\langle g^{G}\rangle)$ divides $o(g)$, the assumption $(o(g),70)=1$ ensures that $o(\overline{g}) \notin \{2,5,7\}$. The result then follows directly from Theorem B.
\end{proof}

\begin{corollary}
Let $G$ be a finite group possessing a conjugacy class $g^{G}$ of size $3p^{n}$ for some prime $p$ and integer $n\ge1$. If $G$  contains no subquotient isomorphic to $A_{5}$, $A_{6}$, $A_{7}$, or $\mathrm{PSL}(2,7)$, then $\langle g^{G}\rangle$ is solvable.
\end{corollary}

\begin{proof}
It is an immediate consequence of Theorem B, as non-solvability would require $G$ to contain one of the specified simple groups as a subquotient.
\end{proof}

\bigskip
\noindent
{\bf Acknowledgements}
This work was supported by the National Natural Science Fund of China (No.  12071181). The second author gratefully acknowledges the warm hospitality received during his research stay at the Universitat Jaume I, Spain, in autumn 2026. 

\bibliographystyle{plain}

\begin{thebibliography}{1}


\bibitem{BS}
A. Beltr\'an and C. Shao,
Finite groups with conjugacy classes of size $2p^n$, arxiv:2610.18193.

\bibitem{BGL}
N. Burgoyne, R. Griess, and R. Lyons,
\emph{Maximal subgroups and automorphisms of Chevalley groups},
Pacific J. Math. \textbf{71} (1977), no. 2, 365--403.

\bibitem{BG}
T. C. Burness and M. Giudici,
\emph{Classical Groups, Derangements and Primes},
Australian Mathematical Society Lecture Series, vol. 25,
Cambridge University Press, Cambridge, 2016.

\bibitem{Camina}
A.~R. Camina and R.~D. Camina, The influence of conjugacy class sizes on the structure of finite groups: a survey, Asian-Eur. J. Math. \textbf{4} (2011), no.~4, 559--588.

\bibitem{Con}
J. H. Conway, R. T. Curtis, S. P. Norton, R. A. Parker, and R. A. Wilson,
\emph{Atlas of Finite Groups},
Oxford University Press, London, 1985.

 \bibitem{Her}
M. Herzog, On finite simple groups of order divisible by three primes only. J. Algebra 10(3), 383–388 (1968).


 \bibitem{Hup67}
B. Huppert,  \textit{ Endliche Gruppen I}. Springer-Verlag, New York-Berlin (1967).

\bibitem{Huppert}
B. Huppert, \textit{Character Theory of Finite Groups}, De Gruyter Expositions in Mathematics, vol. 25, Walter de Gruyter \& Co., Berlin, 1998.

\bibitem{Kaz}
L. S. Kazarin,
\emph{ Burnside's $p^\alpha$-lemma},
 Math. Notes, \textbf{48}, 749–751 (1990).

\bibitem{KL}
P. B. Kleidman and M. W. Liebeck,
\emph{The Subgroup Structure of the Finite Classical Groups},
LMS Lecture Note Series, vol. 129,
Cambridge University Press, Cambridge, 1990.


\bibitem{KS}
H. Kurzweil and B. Stellmacher,
\emph{The Theory of Finite Groups. An Introduction},
Springer-Verlag, Berlin-Heidelberg-New York, 2004.


\bibitem{LL}
C. H. Li and X. H. Li,
\emph{On permutation groups of degree a product of two prime powers},
Comm. Algebra \textbf{42} (2014), no. 11, 4722--4743.

\bibitem{GAP4}
The GAP Group, 
\emph{GAP -- Groups, Algorithms, and Programming, Version 4.13.1}; 
2024, https://www.gap-system.org.
\end{thebibliography}

\end{document}